\documentclass[11pt]{amsart}
\usepackage{amsmath,mathrsfs,cite}

\usepackage[ps2pdf,colorlinks=true,urlcolor=blue,
citecolor=red,linkcolor=blue,linktocpage,pdfpagelabels,bookmarksnumbered,bookmarksopen]{hyperref}
\usepackage[english]{babel}

\usepackage[left=2.7cm,right=2.7cm,top=2.6cm,bottom=2.6cm]{geometry}

\newcommand{\R}{{\mathbb R}}
\newcommand{\T}{{\mathbb T}}
\newcommand{\eps}{\varepsilon}

\numberwithin{equation}{section}
\newtheorem{theorem}{Theorem}[section]

\newtheorem{remark}[theorem]{Remark}
\newtheorem{corollary}[theorem]{Corollary}

\theoremstyle{definition}

\title{On a result by Boccardo-Ferone-Fusco-Orsina}

\author{Marco Squassina}
\address{Dipartimento di Informatica
\newline\indent
Universit\`a degli Studi di Verona
\newline\indent
C\'a Vignal 2, Strada Le Grazie 15
\newline\indent
I-37134 Verona, Italy}
\email{marco.squassina@univr.it}

\thanks{Research supported by PRIN: {\em Metodi Variazionali e Topologici
nello Studio di Fenomeni non Lineari}}

\begin{document}
	

\subjclass[2010]{35A15; 35B06; 74G65; 35B65}

\keywords{Symmetry, non-convex problems, regularity of minimizing sequences.}

\begin{abstract}
Via a symmetric version of Ekeland's principle recently obtained by the author we improve, 
in a ball or an annulus, a result of Boccardo-Ferone-Fusco-Orsina on the properties of minimizing 
sequences of functionals of calculus of variations in the non-convex setting.
\end{abstract}

\maketitle

\section{Introduction}
In the study of non-convex minimization problems \cite{ekeland2} of calculus of variations, 
the idea of selecting minimizing sequences with nice properties to guarantee the convergence 
towards a minimizer can be traced back to Hilbert and Lebesgue \cite{hilbert,lebesgue}.\
In \cite{squassinaJFPTA}, the author has recently obtained an abstract symmetric version of the celebrated
Ekeland's variational principle \cite{ekeland1} for lower semi-continuous 
functionals, which is probably one of the main tools to perform the selection procedure
indicated above. More precisely, the new enhanced Ekeland type principle is able
to select points which are not only almost critical, in a suitable sense,  
but also almost symmetric, provided that the functional does not increase under polarizations \cite{baernst}.
In turn, under rather mild assumptions, starting from a given minimizing sequence one can detect a new
minimizing sequence enriched with very nice features. The 
additional symmetry characteristics play a r\v ole also in non-compact problems, 
providing compactifying effects. In 1999, Boccardo-Ferone-Fusco-Orsina \cite{bffo} 
considered functionals $J:W^{1,p}_0(\Omega)\to\R$ of calculus of variations,
$$
J(u)=\int_\Omega j(x,u,Du),\quad\,\, u\in W^{1,p}_0(\Omega),
$$
with no convexity assumption on $\xi\mapsto j(x,s,\xi)$ and showed that, by merely relying upon some classical \cite{ladura} growth
estimates on the integrand $j(x,s,\xi)$, the existence of minimizing sequences with enhanced 
smoothness can be obtained by combining the application of the classical Ekeland's principle with 
a priori estimates (cf.\ \cite[Lemmas 2.3 and 2.6]{bffo}) based upon suitable Gehring-type lemmas \cite{giaqu}. 
We also refer the reader to \cite{marcsbord} for other results in the same spirit.

The main goal of the present note is to highlight that, if we restrict 
the attention to the case where $\Omega$ is either a ball or an annulus of $\R^N$ and $J$
decreases upon polarizations, then arguing as in \cite{bffo} but using the Ekeland's principle from 
\cite{squassinaJFPTA}, even more special minimizing sequences can be detected. More precisely,
consider $1<p<N$, let $p^*$ denote the critical Sobolev exponent and let 
$j:\Omega\times\R\times\R^N\to\R$ be a Carath\'eodory function such that
\begin{equation}
	\label{growth}
\alpha |\xi|^p-\varphi_2|s|^{\gamma_2}\leq j(x,s,\xi)\leq \beta |\xi|^p+\varphi_0+\varphi_1|s|^{\gamma_1}, 
\end{equation}
for a.e. $x\in \Omega$ and every $(s,\xi)\in \R\times\R^N$, for some $\alpha,\beta>0$,
$$
\varphi_0\in L^{r_0}(\Omega),\,\,\, r_0>1,  
\quad\,\,
\varphi_1\in L^{r_1}(\Omega),\,\,\, r_1>N/p,    
\quad\,\, 
\varphi_2\in L^{r_2}(\Omega),\,\,\, r_2>N, 
$$
and $\gamma_1,\gamma_2$ satisfying
$$
 0\leq \gamma_1<p^*\frac{r_1-1}{r_1},\qquad 
0\leq \gamma_2<\min\Big\{p,\frac{N}{N-1}\frac{r_2-1}{r_2}\Big\}.
$$
We consider the following classes of half-spaces in $\R^N$
\begin{align*}
	& {\mathcal H}_*:=\{\text{$H\subset\R^N$ is a half-space with $0\in H$}\},\quad  \text{if $\Omega$ is a ball},  \\
	& {\mathcal H}_*:=\{\text{$H\subset\R^N$ is a half-space with $\R^+\times\{0\}\subset H$ and
	 $0\in \partial H$}\},\quad \text{if $\Omega$ is an annulus},
\end{align*}
For any nonnegative measurable function $u$ we define $u^H$ to be the polarization of $u$
with respect a half-space $H\in {\mathcal H}_*$. Moreover, we denote by $u^*$ the Schwarz 
symmetrization (resp.\ the spherical cap symmetrization) if $\Omega$ is a ball (resp.\ if 
$\Omega$ is an annulus). For definitions and properties of these notions, we refer  
to \cite{baernst} and to the references therein.
\vskip3pt
\noindent
In this framework, merely under assumption \eqref{growth}, we have the following

\begin{theorem}
	\label{main}
Assume that $\Omega$ is either a ball or an annulus in $\R^N$ with $N\geq 2$ and
\begin{equation}
	\label{polar-ass}
J(u^H)\leq J(u)\qquad\text{for all $u\in W^{1,p}_{0+}(\Omega)$ and any $H\in{\mathcal H}_*$.}
\end{equation}
Then for an arbitrary minimizing sequence $(u_h)\subset W^{1,p}_{0+}(\Omega)$ for $J$ there exist $q>p$, a new minimizing sequence	
$(v_h)\subset W^{1,p}_0(\Omega)$ for $J$ and continuous mappings $\T_h:W^{1,1}_{0+}(\Omega)\to W^{1,1}_{0+}(\Omega)$ such that 
$\T_h z$ is built from $z$ via iterated polarizations by half-spaces in ${\mathcal H}_*$, such that
$$
\sup_{h\geq 1}\|v_h\|_{W^{1,q}_0(\Omega)}<+\infty,\quad\text{if $r_0<\frac{N}{p}$},\qquad
\sup_{h\geq 1}\|v_h\|_{L^\infty(\Omega)}<+\infty,\quad\text{if $r_0>\frac{N}{p}$},
$$
and, in addition,
$$
\lim_h\|v_h-|v_h|^*\|_{L^{\frac{N}{N-1}}(\Omega)}=0,\quad\,\,\,
\limsup_h\|v_h-u_h\|_{W^{1,1}_0(\Omega)}\leq \limsup_h \|\T_h u_h-u_h\|_{W^{1,1}_0(\Omega)}.
$$
\end{theorem}

\noindent
We stress that, under \eqref{growth}, $J$ is bounded from below but, 
since we are not assuming the convexity of $\xi\mapsto j(x,s,\xi)$,
we can by no means conclude that $J$ has a minimum point. Nevertheless, a smooth minimizing
sequence made by almost Schwarz symmetric (for the ball) or almost spherical cap symmetric (for the annulus) points can be constructed.
As it can be readily checked by direct computation, a class of integrands which 
satisfy \eqref{polar-ass} (with the equality in place of the inequality) is, for instance, $j(x,s,\xi)=j_0(s,|\xi|)$, for some 
continuous function $j_0:\R\times\R^+\to\R$. 
Observe also that, of course, from the last conclusion of Theorem \ref{main}, the limit $v$ of
$(v_h)$ must be Schwarz symmetric, namely $v=v^*$.


\begin{remark}\rm
We conclude with an important remark, which is probably one of the main reasons why 
the conclusion of Theorem~\ref{main} is rather powerful in the non-convex framework.
Should one {\em additionally} assume that $\xi\mapsto j(x,s,\xi)$ is convex, it is then often
the case that a functional which satisfies \eqref{polar-ass}, fulfills in turn the corresponding symmetrization inequality 
$J(u^*)\leq J(u)$. In such a case, starting from a given minimizing sequence 
$(u_h)\subset W^{1,p}_{0+}(\Omega)$ for $J$ one has that $(u_h^*)\subset W^{1,p}_{0+}(\Omega)$ is a minimizing sequence too  
and it is then immediate from \cite{bffo} to find a further 
almost symmetric regular minimizing sequence $(v_h)$. 
On the other hand, without the convexity of $j(x,s,\xi)$ in the gradient, to the author knowledge, no symmetrization
inequality is available in the current literature. In some sense, while $J(u^H)\leq J(u)$ is often an algebraic fact, 
$J(u^*)\leq J(u)$ is rather a more geometrical fact. 
\end{remark}

\section{Symmetric Ekeland's principle}
Let $X$ and $V$ be two Banach spaces and $S\subseteq X$. We shall consider two maps $*:S\to V$, $u\mapsto u^*$, 
the symmetrization map, and $h:S\times {\mathcal H}_*\to S$, $(u,H)\mapsto u^H$, the polarization map, ${\mathcal H}_*$ 
being a path-connected topological space. As in \cite{squassinaJFPTA}, we assume the following:
\begin{enumerate}
 \item $X$ is continuously embedded in $V$;
 \item $h$ is a continuous mapping;
\item for each $u\in S$ and $H\in {\mathcal H}_*$ it holds $(u^*)^H=(u^H)^*=u^*$ and $u^{HH}=u^H$;
\item there exists $(H_m)\subset {\mathcal H}_*$ such that, for $u\in S$, $u^{H_1\cdots H_m}$ converges
to $u^*$ in $V$;
\item for every $u,v\in S$ and $H\in {\mathcal H}_*$ it holds
$\|u^H-v^H\|_V\leq \|u-v\|_V$.
\end{enumerate}
Moreover, the mappings $*:S\to V$ and $h:S\times {\mathcal H}_*\to S$ can be extended to $*:X\to V$ 
and $h:X\times {\mathcal H}_*\to S$ respectively by setting 
$u^*:=(\Theta(u))^*$ and $u^H:=(\Theta(u))^H$ for all $u\in X$, where $\Theta:(X,\|\cdot\|_V)\to (S,\|\cdot\|_V)$ is
Lipschitz of constant $C_\Theta$ and such that $\Theta|_{S}={\rm Id}|_{S}$. 
In the above framework, we recall the result from \cite{squassinaJFPTA}.

\begin{theorem}
	\label{mainthm}	
Assume that $f:X\to\R\cup\{+\infty\}$ is a proper and lower 
semi-continuous functional bounded from below such that
\begin{equation}
	\label{assumptionpol}
\text{$f(u^H)\leq f(u)$\,\,\quad for all $u\in S$ and $H\in {\mathcal H}_*$}.
\end{equation}
Let $u\in S$, $\rho>0$ and $\sigma>0$ with
$$
f(u)\leq \inf_X f+\rho\sigma.
$$
Then there exist $v\in X$ and a continuous map $\T_\rho:S\to S$ such that $\T_\rho z$
is built via iterated polarizations of $z$ by half-spaces in ${\mathcal H}_*$ such that
\begin{enumerate}
\item[(a)] $\|v-v^*\|_V< C\rho$; 
\item[(b)] $\|v-u\|_X\leq \rho+\|\T_\rho u-u\|_X;$ 
\item[(c)]   $f(v)\leq f(u)$;   
\item[(d)]   $f(w)\geq f(v)-\sigma \|w-v\|_X,$\quad\text{for all $w\in X,$}
\end{enumerate}
for some positive constant $C$ depending only upon $V,X$ and $\Theta$.
\end{theorem}

\noindent
Let $\Omega$ be either a ball or an annulus of $\R^N$, $N\geq 2$. In particular, by choosing 
$$
X=(W^{1,1}_{0}(\Omega),\|\cdot\|_{W^{1,1}_{0}(\Omega)}),\quad
\|u\|_{W^{1,1}_{0}(\Omega)}=\int_\Omega|Du|,\qquad
S=W^{1,1}_{0+}(\Omega),
$$
as well as
$$
V=(L^{\frac{N}{N-1}}(\Omega),\|\cdot\|_{L^{\frac{N}{N-1}}(\Omega)}),\quad\,\,\, \Theta(u)=|u|,
$$
then (1)-(5) hold true. The following by product, adapted to our purposes, holds true.

\begin{corollary}
	\label{maincor}
	Let $\Omega$ be either a ball or an annulus of $\R^N$ and
	let $J:W^{1,1}_{0}(\Omega)\to\R\cup\{+\infty\}$ be a lower 
	semi-continuous functional bounded from below with
	\begin{equation}
		\label{polar-ass-concr}
	J(u^H)\leq J(u)\qquad\text{for all $u\in W^{1,1}_{0+}(\Omega)$ and $H\in{\mathcal H}_*$.}
	\end{equation}
Let $u\in W^{1,1}_{0+}(\Omega)$ and $\eps>0$ be such that
$$
J(u)\leq \inf_{W^{1,1}_{0}(\Omega)} J+\eps.
$$
Then there exist $v\in W^{1,1}_{0}(\Omega)$ and a continuous map $\T_\eps:W^{1,1}_{0+}(\Omega)\to W^{1,1}_{0+}(\Omega)$ 
such that $\T_\eps z$ is built via iterated polarizations of $z$ by half-spaces in ${\mathcal H}_*$ such that $J(v)\leq J(u)$,
\begin{equation}
	\label{prima}
J(w)\geq J(v)-\sqrt{\eps} \|w-v\|_{W^{1,1}_{0}(\Omega)}, \quad\text{for all $w\in W^{1,1}_{0}(\Omega),$}
\end{equation}
and 
\begin{equation}
	\label{seconda}
\|v-|v|^*\|_{L^{\frac{N}{N-1}}(\Omega)}\leq C\sqrt{\eps},\quad\,\,
\|v-u\|_{W^{1,1}_{0}(\Omega)}\leq \sqrt{\eps}+\|\T_{\eps} u-u\|_{W^{1,1}_{0}(\Omega)},
\end{equation}
for some positive constant $C$.
\end{corollary}

\section{Proof of Theorem~\ref{main}}
\noindent
The argument closely follows \cite[proof of Theorem 3.1, p.128]{bffo}, aiming
to apply Corollary~\ref{maincor} in place of the standard Ekeland's variational principle \cite{ekeland1}.
We will denote by $C$ a generic positive constant which may vary from line to line.
Taking into account that $\gamma_2<p$ and $\gamma_2r_2'<N/(N-1)$,
it easily follows that $\tilde J:W^{1,1}_{0}(\Omega)\to \R\cup\{+\infty\}$,
$$
\tilde J(u)=
\begin{cases}
J(u) & \text{if $u\in W^{1,p}_{0}(\Omega)$,}  \\
+\infty  & \text{if $u\in W^{1,1}_{0}(\Omega)\setminus W^{1,p}_{0}(\Omega)$,}
\end{cases}
$$ 
is a lower semi-continuous functional bounded from below. 
In light of assumption \eqref{polar-ass}, we have
\begin{equation*}
\tilde J(u^H)\leq \tilde J(u)\qquad\text{for all $u\in W^{1,1}_{0+}(\Omega)$ and $H\in{\mathcal H}_*$.}
\end{equation*}
Hence, we are in the framework of Corollary~\ref{maincor}.
Given a minimizing sequence $(u_h)\subset W^{1,p}_{0+}(\Omega)$ for $J$, let $(\eps_h)\subset (0,1]$ 
be such that $\eps_h\to 0$ as $h\to\infty$ and 
\begin{equation}
	\label{minimiseq}
J(u_h)\leq \inf_{W^{1,p}_{0}(\Omega)} J+\eps_h,\quad\,\,\, \text{for any $h\geq 1$}.
\end{equation}
Since the infimum of $J$ over $W^{1,p}_{0}(\Omega)$ equals the infimum of $\tilde J$ over $W^{1,1}_{0}(\Omega)$, it holds
\begin{equation*}
\tilde J(u_h)\leq \inf_{W^{1,1}_{0}(\Omega)} \tilde J+\eps_h,\quad\,\,\, \text{for any $h\geq 1$}.
\end{equation*}
Then, by applying Corollary~\ref{maincor} to $\tilde J$, $u_h$ and $\eps_h$, for any $h\geq 1$, there exists $v_h\in W^{1,1}_{0}(\Omega)$ 
such that $J(v_h)=\tilde J(v_h)\leq \tilde J(u_h)=J(u_h)$ and, for any $h\geq 1$,
\begin{equation}
	\label{addconclu}
\|v_h-|v_h|^*\|_{L^{\frac{N}{N-1}}(\Omega)}\leq C\sqrt{\eps_h},\quad\,\,
\|v_h-u_h\|_{W^{1,1}_{0}(\Omega)}\leq \sqrt{\eps_h}+\|\T_{\eps_h} u_h-u_h\|_{W^{1,1}_{0}(\Omega)},
\end{equation}
for some continuous maps $\T_{\eps_h}:W^{1,1}_{0+}(\Omega)\to W^{1,1}_{0+}(\Omega)$ as well as, for any $h\geq 1$,
\begin{equation*}
\tilde J(v_h)\leq \tilde J(w)+\sqrt{\eps_h} \int_\Omega |Dw-Dv_h|, \quad\text{for all $w\in W^{1,1}_{0}(\Omega)$,}
\end{equation*}
that is, being $\tilde J(v_h)<+\infty$,
\begin{equation}
	\label{conclquasi}
\int_\Omega j(x,v_h,Dv_h)\leq \int_\Omega j(x,w,Dw)+\sqrt{\eps_h} \int_\Omega |Dw-Dv_h|, \quad\text{for all $w\in W^{1,p}_{0}(\Omega)$.}
\end{equation}
Observe that $(u_h)$ is bounded in $W^{1,p}_0(\Omega)$ since \eqref{minimiseq} and \eqref{growth} yield
\begin{equation}
	\label{boundestimm}
\alpha\|Du_h\|^p_{L^p(\Omega)}\leq C+C\|\varphi_2\|_{L^{r_2}(\Omega)}\|Du_h\|_{L^p(\Omega)}^{\gamma_2}, \quad\,\,  (\gamma_2<p).
\end{equation}
In turn, $(v_h)$ is bounded in $W^{1,1}_{0}(\Omega)$, since 
by the second inequality of \eqref{addconclu}, it holds
\begin{align*}
\|v_h\|_{W^{1,1}_{0}(\Omega)} &\leq \|v_h-u_h\|_{W^{1,1}_{0}(\Omega)}+\|u_h\|_{W^{1,1}_{0}(\Omega)} \\
& \leq\sqrt{\eps_h}+\|\T_{\eps_h} u_h-u_h\|_{W^{1,1}_{0}(\Omega)}+\|u_h\|_{W^{1,1}_{0}(\Omega)} \\
& \leq\sqrt{\eps_h}+3\|u_h\|_{W^{1,1}_{0}(\Omega)}\leq 1+C\|u_h\|_{W^{1,p}_{0}(\Omega)}\leq C.
\end{align*}
In the last line, we exploited the fact that, by construction of $\T_{\eps_h}$, for any $h\geq 1$,
$$
\|\T_{\eps_h} u_h\|_{W^{1,1}_{0}(\Omega)}=\int_\Omega \big|Du_h^{H_0\cdots H_{m_{\eps_h}}}\big|
=\int_\Omega \big|Du_h^{H_0\cdots H_{m_{\eps_h}-1}}\big|=\dots=\int_\Omega |Du_h|.
$$
In conclusion, $(v_h)$ is bounded in $W^{1,p}_0(\Omega)$ since by $J(v_h)\leq C$, \eqref{growth} and $\gamma_2r_2'<N/(N-1)$
\begin{equation}
	\label{boundestimm-bis}
\alpha\|Dv_h\|^p_{L^p(\Omega)}\leq C+C\|\varphi_2\|_{L^{r_2}(\Omega)}\|v_h\|_{W^{1,1}_{0}(\Omega)}^{\gamma_2}\leq C, 
\end{equation}
and the variational inequality \eqref{conclquasi}, choosing $w=v_h+\varphi$ for a $\varphi\in W^{1,p}_{0}(\Omega)$, yields 
\begin{equation}
	\label{finalsupt}
\int_{{\rm supt}(\varphi)} j(x,v_h,Dv_h)\leq \int_{{\rm supt}(\varphi)} j(x,v_h+\varphi,Dv_h+D\varphi)
+\sqrt{\eps_h} \int_{{\rm supt}(\varphi)}|D\varphi|, 
\end{equation}
for all $h\geq 1$ and any $\varphi\in W^{1,p}_{0}(\Omega)$.
Once these facts hold, the boundedness of $(v_h)$ in $W^{1,q}_0(\Omega)$ (case $r_0<N/p$) for some $q>p$ follows
as in \cite[proof of Theorem 3.4]{bffo} using \eqref{growth} in \eqref{finalsupt}. 
The boundedness in $L^\infty(\Omega)$ (case $r_0>N/p$), follows by choosing
$w=\max(-k,\min(k,v_h))\in W^{1,p}_0(\Omega)$ in \eqref{conclquasi} and arguing as in 
\cite[proof of Theorems 3.5]{bffo}. Recalling \eqref{addconclu}, the proof is complete.

\vskip25pt
\noindent
{\bf Acknowledgment.} The note is dedicated to the memory of my beloved mother Maria Grazia.

\bigskip
\medskip


\vskip28pt

\begin{thebibliography}{99}

\bibitem{baernst}
{\sc A.\ Baernstein II}, 
A unified approach to symmetrization, 
Sympos. Math., Cambridge U Press (1994), 47--91.
	
\bibitem{bffo}
{\sc L.\ Boccardo, V.\ Ferone, N.\ Fusco, L.\ Orsina}, 
Regularity of minimizing sequences for functionals 
of the calculus of variations via the Ekeland principle, 
{\em Diff. Integral Equations} {\bf 12} (1999), 119--135.

\bibitem{ekeland1}
{\sc I.\ Ekeland}, 
On the variational principle, 
{\em J. Math. Anal. Appl.} {\bf 47} (1974), 324--353.

\bibitem{ekeland2}
{\sc I.\ Ekeland}, 
Nonconvex minimization problems, 
{\em Bull. Amer. Math. Soc.} {\bf 1} (1979), 443--474.

\bibitem{giaqu}
{\sc M.\ Giaquinta}, 
Multiple integrals in the calculus of variations and nonlinear elliptic systems. 
Annals of Mathematics Studies, {\bf 105}. Princeton University Press, Princeton, NJ, 1983. 

\bibitem{hilbert}
{\sc D.\ Hilbert},
Uber das dirichletsche prinzip, 
{\em Jber. Deut. Math. Ver.} {\bf 8} (1900), 184--188. 

\bibitem{ladura}
{\sc O.A.\ Ladyzhenskaya, N.N.\ Ural'ceva}, 
Linear and quasi-linear elliptic equations, New York 1968.

\bibitem{lebesgue}
{\sc H.\ Lebesgue},
Sur le probl\'ems de Dirichlet, 
{\em Rend. Circ. Mat. Palermo} {\bf 24} (1907), 371--402.

\bibitem{marcsbord}
{\sc P.\ Marcellini, C.\ Sbordone}, 
Semicontinuity problems in the calculus of variations, 
{\em Nonlinear Anal.} {\bf 4} (1980), 241--257.

\bibitem{squassinaJFPTA}
{\sc M.\ Squassina}, 
On Ekeland's variational principle, 
{\em J. Fixed Point Theory Appl.} (2011), to appear.

\end{thebibliography}
\end{document}